\documentclass[12pt]{article}

\usepackage{hyperref}
\usepackage{amsmath}
\usepackage{amssymb}
\usepackage{amsfonts}
\usepackage{amsopn}
\usepackage{amsthm}
\usepackage[all]{xy}

\theoremstyle{plain}

\newtheorem{thm}{Theorem}
\newtheorem{cor}{Corollary}

\theoremstyle{definition}
\newtheorem{rem}{Remark}

\newcommand{\zz}{\mathbb{Z}}   

\DeclareMathOperator{\CH}{\mathrm{CH}}      
\DeclareMathOperator{\Ch}{\mathrm{Ch}}

\DeclareMathOperator{\im}{\mathrm{Im}}      
\DeclareMathOperator{\cd}{\mathrm{cd}}  

\title{Canonical $p$-dimensions of algebraic groups and
degrees of basic polynomial invariants}

\author{K.~Zainoulline}

\date{}

\begin{document}

\maketitle

\begin{abstract}
In the present notes we provide a new uniform way to compute 
a canonical $p$-dimension of a split algebraic group $G$ for
a torsion prime $p$
using degrees of basic polynomial invariants described by V.~Kac.
As an application, we compute the canonical $p$-dimensions
for all exceptional simple algebraic groups.
\end{abstract}

The notion of a canonical dimension of an algebraic structure
was introduced by Berhuy and Reichstein \cite{BR}.
For a split algebraic group $G$ and its torsion prime $p$ 
the canonical $p$-dimension of $G$ was studied by 
Karpenko and Merkurjev in \cite{KM}.
In particular, this invariant was shown to be related with
the size of the image of the characteristic map 
\begin{equation}
\phi_G: S^*(\hat{T}) \to \CH^*(X),
\end{equation}
where $\hat{T}$ is the character group of a maximal split torus $T$,
$X$ is the variety of complete flags and $S^*$ stands for the
 symmetric algebra.
Namely, one has the following formula for 
the canonical $p$-dimension of a group $G$
\begin{equation}
\cd_p(G)=\min \{i \mid \overline{\Ch_i}(X)\neq 0\},
\end{equation}
where $\overline{\Ch_i}(X)$ stands for the image of $\phi_G$ in the modulo $p$
Chow group $\Ch(X)=\CH(X)/p\cdot \CH(X)$ (see \cite[Theorem~6.9]{KM}).
Using this remarkable fact together with the explicit description
of the image $R_p=\overline{\Ch}(X)$
Karpenko and Merkurjev computed the canonical $p$-dimensions
for all classical algebraic groups $G$ (see \cite[Section~8]{KM}).

The goal of the present notes is to relate
the work by V.~Kac \cite{Kc}
devoted to the study of the $p$-torsion part of 
the cohomology ring of an algebraic group $G$
with the canonical $p$-dimensions of $G$, hence, providing
a different and uniform approach
of computing those invariants. As a consequence, we compute
canonical $p$-dimensions for all exceptional
algebraic groups, hence, completing the computations started in \cite{KM}.

\begin{thm}\label{main} Let $G$ be a split simple algebraic group of rank $n$
and $p$ be an odd torsion prime. 
Then
$$
\cd_p(G)=N + n - (d_{1,p}+ d_{2,p}+\ldots+d_{n,p}),
$$
where $N$ stands for the number of positive roots of $G$ and 
integers $d_{1,p},\ldots,d_{n,p}$ are 
the degrees of basic polynomial invariants modulo $p$.
\end{thm}

\begin{proof} Consider the characteristic map (1) modulo $p$
\begin{equation}
(\phi_G)_p: S^*(\hat{T})\otimes_\zz \zz/p\zz \to \Ch^*(X)
\end{equation}
According to \cite[Cor.~3.9.(ii)]{NNS} 
the kernel of this map $I_p$ is generated by a regular sequence
of $n$ homogeneous polynomials of degrees $d_{1,p},\ldots,d_{n,p}$.

Recall that a Poincare polynomial $P(A,t)$ for a graded module $M^*$ 
over a field $k$ is defined to be
$\sum_i \dim_k M^i \cdot t^i$ (see \cite{St}).
Hence, the Poincare polynomial for the $\zz/p\zz$-module 
$R_p$ is equal to
\begin{equation}
P(R_p,t)= \prod_{i=1}^n \frac{1-t^{d_{i,p}}}{1-t}.
\end{equation}
Indeed, we identify $R_p$ with the quotient of the polynomial ring
in $n$ variables $\zz/p\zz[\omega_1,\ldots,\omega_n]$ modulo the ideal $I_p$.
The formula (4) then follows immediately by \cite[Cor.~3.3]{St}.

According to (2) the canonical $p$-dimension is equal to the difference
$\dim(X)- \deg P(R_p,t)=N - \sum_{i=1}^n (d_{i,p}-1)$.
\end{proof}

\setcounter{rem}{\value{thm}}\addtocounter{thm}{1}
\begin{rem} The fact that the ideal
$I_p$ is generated by a regular sequence of elements 
was extensively  used in \cite{Kc}. Unfortunately, 
the original proof of it (provided in \cite{Kc}) contains a mistake.
Another proof of this fact which works only for odd torsion primes
can be found in \cite{NNS}.
\end{rem}

The next theorem relates $\cd_p(G)$ with the Chow group of $G$ modulo $p$.

\begin{thm}\label{group} Let $G$ be a split simple group and $p$ be its torsion prime.
Then
$$
\cd_p(G)= \max \{i \mid \Ch^i(G)\neq 0 \}
$$
\end{thm}

\begin{proof}
It is known that $\Ch(G)=\Ch(X)/J_p$ (see \cite{Gr}),
where $J_p$ is the ideal generated by the non-constant part of $R_p$.
Since $\Ch(X)$ is a free $R_p$-module (see \cite[Appendix]{Kc}), we have
$P(\Ch(X)/J_p,t)\cdot P(R_p,t)=P(\Ch(X),t)$ and, hence, for the degrees
$\deg(P(\Ch(G),t))+\deg(P(R_p,t))=\dim(X)$. The proof is completed since
$\deg(P(R_p,t))=\dim(X)-\cd_p(G)$.
\end{proof}

\setcounter{cor}{\value{thm}}\addtocounter{thm}{1}
\begin{cor}\label{pexc}
Let $d_1,d_2,\ldots,d_r$ be the set of $p$-exceptional degrees 
(introduced in \cite[Thm.~3]{Kc})
for a group $G$ and its odd torsion prime $p$.
Then the canonical
$p$-dimension of $G$ is equal to the sum
$$
\cd_p(G) = \sum_{i=1}^r d_i'\cdot (p^{k_i}-1),
$$
where the integers $d_i'$ and $p^{k_i}$ are the factors of the decompositions
$d_i=d_i'\cdot p^{k_i}$, $p \nmid d_i'$.
\end{cor}

\begin{proof} Follows by the last isomorphism of \cite[Theorem 3.(ii)]{Kc}.
\end{proof}

\begin{cor} We obtain the following values for the canonical 
$p$-dimensions
of groups of types $F_4$, $E_6$, $E_7$ and $E_8$
(here $G^{sc}$ and $G^{ad}$ stand for
the simply-connected and adjoint forms of a group $G$)

\begin{tabular}{ccc}
$\cd_2 F_4 = 3$, & $\cd_3 F_4 = 8$ & \\
$\cd_2 E_6 = 3$, & $\cd_3 E_6^{sc} = 8$, & $\cd_3 E_6^{ad} = 16$ \\
$\cd_2 E_7^{sc} = 17$, & $\cd_2 E_7^{ad} = 18$ & $\cd_3 E_7 = 8$ \\
$\cd_2 E_8 = 60$, & $\cd_3 E_8= 28$ & $\cd_5 E_8= 24$ 
\end{tabular}
\end{cor}

\begin{proof}
The case of odd torsion primes follows immediately by Theorem~\ref{main} or
Corollary~\ref{pexc} and Table~2 of \cite{Kc}.

The case $p=2$ can be computed using Theorem~\ref{group} as follows.
Consider the canonical map $\pi: G \to X$. 
Note that the Chow group $\Ch(G)$ can be identified with the
image $\im \pi^*$ of the induced pull-back (see \cite{Gr}). 
According to \cite[Thm.~1.1 and Lem.~1.3]{IKT}, 
the cohomology ring of $G$ can be represented as the tensor
product of two algebras
\begin{equation}\label{cohom}
H(G;\zz/p\zz)\cong \im \pi^* \otimes \Delta(a_1,\ldots, a_n),
\end{equation}
where each $a_i$ is of odd degree, all elements of $\im \pi^*$ are of even degree and
$\Delta(a_1,\ldots,a_n)$ denotes the submodule
spanned by the simple monomials $a_1^{\epsilon_1},\ldots,a_n^{\epsilon_n}$
($\epsilon_i=0$ or $1$) which are linearly independent.
Knowing this representation and the cohomology ring of $G$, one
immediately obtains 
$$
cd_p(G)=\tfrac{1}{2}(\deg P(H(G;\zz/p\zz),t) - \sum_{i=1}^n \deg(a_i)).
$$
To finish the proof 
observe that the cohomology ring of an exceptional algebraic
group modulo 2 and degrees of the elements $a_i$ 
can be found in the literature (see the references of paper \cite{IKT}).
\end{proof}

{\it Acknowledgements.} 
I am grateful to Burt Totaro and Larry~Smith for very useful 
comments concerning paper \cite{Kc}.

\bibliographystyle{chicago}

\begin{thebibliography}{99}



\bibitem{BR} Berhuy, G., Reichstein, Z. On the notion of canonical dimension
for algebraic groups. {\it Advances in Math.}, in press, doi: 10.1016/j.aim.2004.12.004.

\bibitem{Gr} Grothendieck, A., Torsion homologique et sections rationnelles in 
{\it Anneaux de Chow et applications}, S\'eminaire C.~Chevalley 2, 1958.

\bibitem{IKT} Ishitoya, K., Kono, A., Toda, H.
Hopf Algebra Structure of mod 2 Cohomology of Simple Lie Groups.
{\it Publ. RIMS, Kyoto Univ.} 12 (1976), 141--167.

\bibitem{Kc} Kac, V. Torsion in cohomology of compact Lie groups and
Chow rings of reductive algebraic groups. {\it Invent. Math.} 80 (1985), 69--79.

\bibitem{KM} Karpenko, N., Merkurjev, A. Canonical $p$-dimension of algebraic groups.
{\it Advances in Math.}, in press, doi: 10.1016/j.aim.2005.07.013.

\bibitem{St} Stanley, R. Hilbert functions of graded algebras. {\it Advances in Math.} 28 (1978), 57--83.

\bibitem{NNS} Neumann, F., Neusel, M., Smith, L. 
Rings of generalized and stable
invariants and classifying spaces of compact Lie groups in {\it Higher homotopy structures in top. and math. physics (Poughkeepsie, NY, 1996)}, Contemp. Math., 227 (1999), 267--285.

\end{thebibliography}

\end{document}